\documentclass[a4paper,final,oneside]{article}%
\usepackage{amsmath}
\usepackage{graphicx}
\usepackage{amsfonts}
\usepackage{amssymb}%
\begin{document}

\title{The near--Hopf ring structure on the integral cohomology of 1--connected Lie groups}
\author{Haibao Duan\\\noindent Institute of Mathematics, Chinese Academy of Sciences\\dhb@math.ac.cn}
\date{}
\maketitle

\begin{abstract}
Let $G$ be an exceptional Lie group with a maximal torus $T$. Based on
Schubert calculus on the flag manifold $G/T$ we have described the integral
cohomology ring $H^{\ast}(G)$ by explicitely constructed generators in
\cite{[DZ2]}, and determined the structure of $H^{\ast}(G;\mathbb{F}_{p})$ as
a Hopf algebra over the Steenrod algebra in \cite{[DZ3]}. Combining these
works we determine the near--Hopf ring structure on $H^{\ast}(G)$.

\begin{description}
\item \textsl{2000 Mathematical Subject Classification:} 55T10

\item \textsl{Key words and phrases:} Lie groups; cohomology; Hopf algebra

\end{description}
\end{abstract}

\section{Introduction}

Let $G$ be a compact $1$--connected Lie group with multiplication $\mu:G\times
G\rightarrow G$ and integral cohomology ring $H^{\ast}(G)$. The induced ring
map\textsl{ }

\begin{enumerate}
\item[(1.1)] $\mu^{\ast}:H^{\ast}(G)\rightarrow H^{\ast}(G\times G)$
\end{enumerate}

\noindent will be referred to as the\textsl{ near--Hopf ring structure} on
$H^{\ast}(G)$. In \cite[Theorem 6]{[DZ2]} we have described the rings
$H^{\ast}(G)$ for all exceptional Lie groups $G$ by explicitly constructed
\textsl{primary generators}. In this sequel to \cite{[DZ2]} we determine the
action of $\mu^{\ast}$ on $H^{\ast}(G)$ with respect to these generators.

We shall start with a brief introduction in \S 2 on two sets of generators for
$H^{\ast}(G)$; $\varrho_{i_{1}},\cdots,\varrho_{i_{n}}$ and $\mathcal{C}_{I}$
with $I$ a certain multi--index, constructed explicitly in \cite{[DZ2]}. The
elements $\varrho_{i}$ have infinite order and their square free products form
a $\mathbb{Z}$--basis for the free part of $H^{\ast}(G)$. The $\mathcal{C}%
_{I}$ come from the Bockstein of a distinguished set of elements in $H^{\ast
}(G;\mathbb{F}_{p})$ with $p=2,3$ or $5$.

For a prime $p$ the multiplication $\mu$ induces also the \textsl{co--product}

\begin{enumerate}
\item[(1.2)] $\mu_{p}^{\ast}:H^{\ast}(G;\mathbb{F}_{p})\rightarrow H^{\ast
}(G;\mathbb{F}_{p})\otimes H^{\ast}(G;\mathbb{F}_{p})$
\end{enumerate}

\noindent on the algebra $H^{\ast}(G;\mathbb{F}_{p})$ by virtue of the Kunneth
formula. It furnishes $H^{\ast}(G;\mathbb{F}_{p})$ with the structure of a
\textsl{Hopf algebra}. In \S 3 we establish a "\textsl{pull--back formula}"
that reduces the computation of $\mu^{\ast}(z)$ to that of $\mu_{p}^{\ast
}(z\operatorname{mod}p)$ for a generators $z=\varrho_{i}$ or $\mathcal{C}_{I}%
$, and deduce from \cite[Theorem 2; Theorem 4.1]{[DZ3]} a presentation of the
Hopf algebra $H^{\ast}(G;\mathbb{F}_{p})$. 

With these preparation the near--Hopf ring $H^{\ast}(G)$ for an exceptional
$G$ is obtained and presented in Theorems 1--5 in \S 4.

The author feels grateful to M. Mimura for encouraging him to work out these
results. Earlier in the 1950's Borel calculated the integral cohomology ring
$H^{\ast}(G)$ for $G=G_{2},F_{4}$ \cite{[B1], [B2]}. However, the action of
$\mu^{\ast}$ was absent even in these cases. On the other hand, combining our
results with previous knowledge for the classical groups $G=SU(n)$, $Sp(n)$,
$Spin(n)$ by Borel and Pitties \cite{[B2], [P]} completes the project to
determine the near--Hopf ring structure of all compact $1$--connected Lie groups.

\section{The generators for the ring $H^{\ast}(G)$}

Assume that $G$ is a compact $1$--connected Lie group with a maximal torus
$T\subset G$, and let $\mathbb{F}$ be either the ring $\mathbb{Z}$ of
integers, the field $\mathbb{Q}$ of rationals, or one of the finite fields
$\mathbb{F}_{p}$. In view of the Leray--Serre spectral sequence $\{E_{r}%
^{\ast,\ast}(G;\mathbb{F}),d_{r}\}$ of the fibration

\begin{enumerate}
\item[(2.1)] $T\rightarrow G\overset{\pi}{\rightarrow}G/T$
\end{enumerate}

\noindent with coefficients in $\mathbb{F}$, we recall from \cite{[DZ2]} the
two sets $\varrho_{i}$'s and $\mathcal{C}_{I}$'s of generators for $H^{\ast
}(G)$ mentioned in \S 1.

The next result was due to Leray \cite{[L]} for $\mathbb{F=Q}$, extended to
$\mathbb{F=F}_{p}$ by Serre \cite{[S]}, conjectured for $\mathbb{F=Z}$ by
Ka\v{c} \cite{[K]} and Marlin \cite{[M]}.

\bigskip

\noindent\textbf{Lemma 2.1 }(\cite[Remark 5.3]{[DZ2]})\textbf{. }As graded
groups $H^{\ast}(G;\mathbb{F})=E_{3}^{\ast,\ast}(G;\mathbb{F})$.$\square$

\bigskip

Lemma 2.1 grants us with the ready--made decomposition

\begin{enumerate}
\item[(2.2)] $H^{\ast}(G;\mathbb{F})=E_{3}^{\ast,0}(G;\mathbb{F})\oplus
E_{3}^{\ast,1}(G;\mathbb{F})\oplus\cdots\oplus E_{3}^{\ast,n}(G;\mathbb{F})$.
\end{enumerate}

\noindent Concrete presentations for the lower terms $E_{3}^{\ast
,0}(G;\mathbb{Z})$ and $E_{3}^{\ast,1}(G;\mathbb{Z})$ bring us certain
elements in $H^{\ast}(G)$ that suffice to generate the ring $H^{\ast}(G)$. We
start by introducing for each exceptional $G$ a set of \textsl{special
Schubert classes} $y_{i}\in H^{i}(G/T)$ on $G/T$ by their\textsl{ Weyl
coordinates} \cite{[DZ1]} in the table below:

\begin{enumerate}
\item[(2.3)]
\begin{tabular}
[c]{l|l|l|l}\hline\hline
$G/T$ & $G_{2}/T$ & $F_{4}/T$ & $E_{n}/T,\text{ }n=6,7,8$\\\hline
$y_{6}$ & $\sigma_{\lbrack1,2,1]}$ & $\sigma_{\lbrack3,2,1]}$ & $\sigma
_{\lbrack5,4,2]}\text{, }n=6,7,8$\\
$y_{8}$ &  & $\sigma_{\lbrack4,3,2,1]}$ & $\sigma_{\lbrack6,5,4,2]}\text{,
}n=6,7,8$\\
$y_{10}$ &  &  & $\sigma_{\lbrack{7,6,5,4,2}]}\text{, }n=7,8$\\
$y_{12}$ &  &  & $\sigma_{\lbrack1,3,6,5,4,2]}\text{, }n=8$\\
$y_{18}$ &  &  & $\sigma_{\lbrack1,5,4,3,7,6,5,4,2]}$,$\text{ }n=7,8$\\
$y_{20}$ &  &  & $\sigma_{\lbrack{1,6,5,4,3,7,6,5,4,2}]}\text{, }n=8$\\
$y_{30}$ &  &  & $\sigma_{\lbrack5,4,2,3,1,6,5,4,3,8,7,6,5,4,2]}\text{, }%
n=8$\\\hline\hline
\end{tabular}

\end{enumerate}

\noindent and set $x_{i}:=\pi^{\ast}(y_{i})\in H^{\ast}(G)$. The next result
was shown in \cite[Lemma 3.1]{[DZ2]}.

\bigskip

\noindent\textbf{Lemma 2.2.} The subring $E_{3}^{\ast,0}(G;\mathbb{Z})\subset
H^{\ast}(G)$ has the presentation

\begin{quote}
$E_{3}^{\ast,0}(G_{2};\mathbb{Z})=\mathbb{Z}[x_{6}]/\left\langle 2x_{6}%
,x_{6}^{2}\right\rangle $;

$E_{3}^{\ast,0}(F_{4};\mathbb{Z})=\mathbb{Z}[x_{6},x_{8}]/\left\langle
2x_{6},x_{6}^{2},3x_{8},x_{8}^{3}\right\rangle $;

$E_{3}^{\ast,0}(E_{6};\mathbb{Z})=\mathbb{Z}[x_{6},x_{8}]/\left\langle
2x_{6},x_{6}^{2},3x_{8},x_{8}^{3}\right\rangle $;

$E_{3}^{\ast,0}(E_{7};\mathbb{Z})=\mathbb{Z}[x_{6},x_{8},x_{10},x_{18}%
]/\left\langle 2x_{6},3x_{8},2x_{10},2x_{18},x_{6}^{2},x_{8}^{3},x_{10}%
^{2},x_{18}^{2}\right\rangle $;

$E_{3}^{\ast,0}(E_{8};\mathbb{Z})=\mathbb{Z}[x_{6},x_{8},x_{10},x_{12}%
,x_{18},x_{20},x_{30}]/$

$\quad\left\langle 2x_{6},3x_{8},2x_{10},5x_{12},2x_{18},3x_{20},2x_{30}%
,x_{6}^{8},x_{8}^{3},x_{10}^{4},x_{12}^{5},x_{18}^{2},x_{20}^{3},x_{30}%
^{2}\right\rangle $.$\square$
\end{quote}

Let $(G,p)$\ be a pair with $G$\ an exceptional Lie group and $H^{\ast}%
(G)$\ containing non--trivial $p$--torsion, and let $r(G,p)$ be the set of the
degrees of the basic \textquotedblleft\textsl{generalized Weyl invariants}%
\textquotedblright\ of $G$\ over $\mathbb{F}_{p}$. Precisely we have (see
\cite{[K]})

\begin{enumerate}
\item[(2.4)]
\begin{tabular}
[c]{ll}\hline
$(G,p)$ & $e(G,p)\subset r(G,p)$\\\hline
$(G_{2},2)$ & $\{2,\underline{3}\}$\\
$(F_{4},2)$ & $\{2,\underline{3},8,12\}$\\
$(E_{6},2)$ & $\{2,\underline{3},5,8,9,12\}$\\
$(E_{7},2)$ & $\{2,\underline{3},\underline{5},8,\underline{9},12,14\}$\\
$(E_{8},2)$ & $\{2,\underline{3},\underline{5},8,\underline{9}%
,12,14,\underline{15}\}$\\
$(F_{4},3)$ & $\{2,\underline{4},6,8\}$\\
$(E_{6},3)$ & $\{2,\underline{4},5,6,8,9\}$\\
$(E_{7},3)$ & $\{2,\underline{4},6,8,10,14,18\}$\\
$(E_{8},3)$ & $\{2,\underline{4},8,\underline{10},14,18,20,24\}$\\
$(E_{8},5)$ & $\{2,\underline{6},8,12,14,18,20,24\}$\\\hline
\end{tabular}
.
\end{enumerate}

\noindent Let $e(G,p)$ be the subset of $r(G,p)$ whose elements are
underlined. In [DZ$_{2}$, (2.10)] we have constructed a subset $\{\zeta
_{2s-1}\in E_{3}^{2s-2,1}(G;\mathbb{F}_{p})\mid s\in r(G,p)\}$, whose elements
are called $p$\textsl{--primary generators on }$H^{\ast}(G;\mathbb{F}_{p})$,
that satisfys the following properties (\cite[Theorem 1]{[DZ2]}).

\bigskip

\noindent\textbf{Lemma 2.3.} \textsl{One has the additive presentation}

\begin{quote}
$H^{\ast}(G;\mathbb{F}_{p})=E_{3}^{\ast,0}(G_{2};\mathbb{F}_{p})\otimes
\Delta(\zeta_{2s-1})_{s\in r(G,p)}$
\end{quote}

\noindent\textsl{on which the Bockstein }$\delta_{p}:H^{r}(G;F_{p})\rightarrow
H^{r+1}(G)$\textsl{ is given by}

\begin{quote}
$\delta_{p}(\zeta_{2s-1})=\left\{
\begin{tabular}
[c]{l}%
$-x_{2s}$ if $s\in e(G,p)$;\\
$0$ if $s\notin e(G,p)$,
\end{tabular}
\ \right.  $
\end{quote}

\noindent\textsl{where }$E_{3}^{\ast,0}(G;\mathbb{F}_{p})=E_{3}^{\ast
,0}(G;\mathbb{Z})\otimes\mathbb{F}_{p}$\textsl{, and where }$\Delta
(\zeta_{2s-1})_{s\in r(G,p)}$\textsl{ is the }$\mathbb{F}_{p}$\textsl{--space
with basis the square free monomials in }$\{\zeta_{2s-1}\mid s\in
r(G,p)\}$\textsl{.}$\square$

\bigskip

\noindent\textbf{Definition 2.4.} For a subset $I\subseteq e(G,p)$ we put
$\mathcal{C}_{I}:=\delta_{p}(\zeta_{I})\in H^{\ast}(G)$, where $\zeta_{I}=%
{\textstyle\prod\limits_{s\in I}}
\zeta_{2s-1}\in E_{3}^{\ast,k}(G;\mathbb{F}_{p})$ with $k$ the cardinality of
$I$.$\square$

\bigskip

For a presentation of $E_{3}^{\ast,1}(G;\mathbb{Z})$ we note that the cup
product in $H^{\ast}(G)$ defines an action of $E_{3}^{\ast,0}(G;\mathbb{Z})$
on $E_{3}^{\ast,1}(G;\mathbb{Z})$. Given a ring $A$ and a finite set
$\{u_{1},\cdots,u_{t}\}$ write $A\{u_{i}\}_{1\leq i\leq t}$ for the free
$A$--module with basis $\{u_{1},\cdots,u_{t}\}$. For an exceptional $G$ with
rank $n$ let $\mathcal{D}_{G}=\{d_{1},\cdots,d_{n}\}$ be the set of degrees of
basic $W$--invariants of $G$ over $\mathbb{Q}$. As is classical we have (e.g.
\cite{[K]})

\begin{enumerate}
\item[(2.5)]
\begin{tabular}
[c]{l|l}\hline\hline
Type of $G$ & $\mathcal{D}_{G}$\\\hline
$G_{2}$ & $\{2,6\}$\\
$F_{4}$ & $\{2,3,8,12\}$\\
$E_{6}$ & $\{2,5,6,8,9,12\}$\\
$E_{7}$ & $\{2,6,8,10,12,14,18\}$\\
$E_{8}$ & $\{2,8,12,14,18,20,24,30\}$\\\hline\hline
\end{tabular}

\end{enumerate}

\noindent The next result is indicated in the proof of \cite[Theorem 2]{[DZ2]}.

\bigskip

\noindent\textbf{Lemma 2.5.} \textsl{For each }$i\in\mathcal{D}_{G}$\textsl{
there is an element }$\varrho_{2i-1}\in E_{3}^{2i-2,1}(G;\mathbb{Z})$\textsl{
of infinite order so that }

\begin{quote}
$E_{3}^{\ast,1}(G;\mathbb{Z})=E_{3}^{\ast,0}(G;\mathbb{Z})\{\varrho
_{2i-1}\}_{i\in\mathcal{D}_{G}}$\textsl{.}
\end{quote}

\noindent\textsl{Moreover, a set of ring generators for }$H^{\ast}(G)$\textsl{
is}

\begin{enumerate}
\item[(2.6)] $\mathcal{G}(G)=\{\varrho_{2i-1},C_{I}\mid i\in\mathcal{D}_{G}%
$\textsl{, }$I\subseteq e(G,p),p=2,3,5\}$\textsl{.}$\square$
\end{enumerate}

\noindent\textbf{Remark 2.6.} In (2.3) the Weyl coordinate of a Schubert class
$y_{i}\in$ $H^{i}(G/T)$ allows one to\textsl{ }construct explicitly the
Schubert variety on $G/T$ Kronecker dual to $y_{i}$ \cite{[DZ1]}.

The subring $E_{3}^{\ast,0}(G;\mathbb{Z})\subset H^{\ast}(G)$ coincides with
the Chow ring of the reductive algebraic group $G^{c}$ corresponding to $G$,
see \cite[Lemma 3.1]{[DZ2]}.$\square$

\section{The pull-back formula}

For convenience we fashion from $\mu^{\ast}$ the \textsl{reduced co--product}

\begin{quote}
$\psi:H^{\ast}(G)\rightarrow H^{\ast}(G\times G)$
\end{quote}

\noindent by $\psi(x)=\mu^{\ast}(x)-(x\otimes1+1\otimes x)$. Its
$\mathbb{F}_{p}$--analogue (resp. $\mathbb{Q}$--analogue) is denoted by
$\psi_{p}$ (resp. $\psi_{0}$). The next result contains a formula that reduces
calculation of $\psi$ to that of $\psi_{p}$.

For a topological space $X$ let $\tau(X)$ be the torsion ideal in the integral
cohomology $H^{\ast}(X)$. For a prime $p$ write $\tau_{p}(X)$ be the
$p$--\textsl{primary component }of $\tau(X)$. The $\operatorname{mod}p$
reduction $H^{r}(X)\rightarrow H^{r}(X;\mathbb{F}_{p})$ is denoted by $r_{p}$.

\bigskip

\noindent\textbf{Lemma 3.1.} \textsl{For an exceptional Lie group }$G$\textsl{
let }$\mathcal{G}(G)$\textsl{ be the set of generators for }$H^{\ast}%
(G)$\textsl{ given in (2.5). Then}

\begin{quote}
\textsl{a) }$\psi(x)\in\tau(G\times G)$\textsl{ for all }$x\in\mathcal{G}%
(G)$\textsl{.}
\end{quote}

\noindent\textsl{Moreover,} \textsl{for }$K=G\times G$\textsl{ or }$G$\textsl{
}

\begin{quote}
\textsl{b) }$\tau(K)=\oplus_{p=2,3,5}\tau_{p}(K)$\textsl{;}

\textsl{c) the reduction }$r_{p}$\textsl{ restricts to an injection }$\tau
_{p}(K)\rightarrow H^{\ast}(K;F_{p})$\textsl{, }
\end{quote}

\noindent\textsl{In particular, for all }$z\in\mathcal{G}(G)$\textsl{ one has}

\begin{enumerate}
\item[(3.1)] $\psi(z)=r_{2}^{-1}\psi_{2}(r_{2}(z))+r_{3}^{-1}\psi_{3}%
(r_{3}(z))+r_{5}^{-1}\psi_{5}(r_{5}(z))$.
\end{enumerate}

\noindent\textbf{Proof.} For $x=$ $\mathcal{C}_{I}$ assertion a) follows from
$x\in\tau(G)$. Assume now that $x=\varrho_{j}\in\mathcal{G}(G)$. The map
$i:H^{\ast}(G)\rightarrow H^{\ast}(G;\mathbb{Q})$ induced by the inclusion
$\mathbb{Z}\rightarrow\mathbb{Q}$ of coefficients clearly preserves the
decomposition (2.2) and therefore, restricts to a homomorphism

\begin{quote}
$i:E_{3}^{\ast,1}(G;\mathbb{Z})=E_{3}^{\ast,0}(G;\mathbb{Z})\{\varrho
_{2i-1}\}_{i\in\mathcal{D}_{G}}\rightarrow E_{3}^{\ast,1}(G;\mathbb{Q})$.
\end{quote}

\noindent Since the subset $E_{3}^{\ast,1}(G;\mathbb{Q})\subset H^{\ast
}(G;\mathbb{Q})$ consists of primitive elements \cite{[L]}, we get
$\psi(\varrho_{j})\in\tau(G\times G)$ from $\psi_{0}(i(\varrho_{j}))=0$. This
shows a).

According to \cite[Theorem 2]{[DZ2]} results in b) and c) hold for any
$1$--connected Lie group $K$. Finally, the commutative diagram

\begin{center}%
\begin{tabular}
[c]{lll}%
$H^{\ast}(G)$ & $\overset{r_{p}}{\rightarrow}$ & $H^{\ast}(G;\mathbb{F}_{p}%
)$\\
$\mu^{\ast}\downarrow$ &  & $\mu_{p}^{\ast}\downarrow$\\
$H^{\ast}(G\times G)$ & $\overset{r_{p}}{\rightarrow}$ & $H^{\ast}(G\times
G;\mathbb{F}_{p})$%
\end{tabular}

\end{center}

\noindent induced by $\mu$ implies that $\operatorname{Im}\psi_{p}\circ
r_{p}=\operatorname{Im}r_{p}\circ\psi\subset H^{\ast}(G\times G;\mathbb{F}%
_{p})$. One obtains (3.1) from a), b) and c), where the validity of
$r_{p}^{-1}$, $p=2,3,5$, follows from the injectivity of $r_{p}:$ $\tau
_{p}(G\times G)\rightarrow H^{\ast}(G\times G;\mathbb{F}_{p})$ by c).$\square$

\bigskip

In order to apply formula (3.1) to determine the action of $\psi$ we need
accounts for

\begin{quote}
i) the reduction $r_{p}:H^{\ast}(G)\rightarrow H^{\ast}(G;\mathbb{F}_{p})$
with respect to the set $\mathcal{G}(G)$ of generators for $H^{\ast}(G)$ and
the presentation of $H^{\ast}(G;\mathbb{F}_{p})$ in Lemma 2.3;

ii) the Hopf algebra structure on $H^{\ast}(G;\mathbb{F}_{p})$ with respect to
its presentation in Lemma 2.3.
\end{quote}

\noindent For i) we have the next result from \cite[Lemma 3.5]{[DZ2]}.

\bigskip

\noindent\textbf{Lemma 3.2.} \textsl{With respect to the set }$\mathcal{G}%
(G)$\textsl{ of generators on }$H^{\ast}(G)$\textsl{ and the set }$\left\{
\zeta_{2s-1}\right\}  _{s\in r(G,p)}$\textsl{ of }$p$\textsl{--primary
generators on }$H^{\ast}(G;\mathbb{F}_{p})$\textsl{, the }$\operatorname{mod}%
p$\textsl{ reduction }$r_{p}:H^{\ast}(G)\rightarrow H^{\ast}(G;\mathbb{F}%
_{p})$\textsl{ is given by}

\begin{enumerate}
\item[(3.2)]
\begin{tabular}
[c]{l|l|l|l}\hline\hline
& $p=2$ & $p=3$ & $p=5$\\\hline
$r_{p}(\varrho_{3})$ & $\zeta_{3}$ & $\zeta_{3}$ & $\zeta_{3}$\\\hline
$r_{p}(\varrho_{9})$ & $\zeta_{9}$ & $\zeta_{9}$ & $\zeta_{9}$\\\hline
$r_{p}(\varrho_{11})$ & $x_{6}\zeta_{5}$ & $-\zeta_{11}$ & $2\zeta_{11}%
$\\\hline
$r_{p}(\varrho_{15})$ & $\zeta_{15}$ & $\zeta_{15}$ & $\zeta_{15}$\\\hline
$r_{p}(\varrho_{17})$ & $\zeta_{17}$ & $\zeta_{17}$ & $\zeta_{17}$\\\hline
$r_{p}(\varrho_{19})$ & $x_{10}\zeta_{9}$ & $-\zeta_{19}$ & $2\zeta_{19}%
$\\\hline
$r_{p}(\varrho_{23})$ & $\zeta_{23}$ & $-x_{8}^{2}\zeta_{7}$ & $2\zeta_{23}%
$\\\hline
$r_{p}(\varrho_{27})$ & $\zeta_{27}$ & $\zeta_{27}$ & $\zeta_{27}$\\\hline
$r_{p}(\varrho_{35})$ & $x_{18}\zeta_{17}$ & $-\zeta_{35}$ & $2\zeta_{35}%
$\\\hline
$r_{p}(\varrho_{39})$ & $x_{10}^{3}\zeta_{9}$ & $-\zeta_{39}$ & $2\zeta_{39}%
$\\\hline
$r_{p}(\varrho_{47})$ & $x_{6}^{7}\zeta_{5}$ & $-\zeta_{47}$ & $2\zeta_{47}%
$\\\hline
$r_{p}(\varrho_{59})$ & $x_{30}\zeta_{29}$ & $-x_{20}^{2}\zeta_{19}$ &
$2x_{12}^{4}\zeta_{11}$\\\hline\hline
\end{tabular}
\textsl{;}

\item[(3.3)] $r_{p}(\mathcal{C}_{I})=\beta_{p}(\zeta_{I})$, $I\subseteq
e(G;p)$,
\end{enumerate}

\noindent\textsl{where }$\beta_{p}=r_{p}\circ\delta_{p}$\textsl{ is the
Bockstein operator on }$H^{\ast}(G;\mathbb{F}_{p})$\textsl{.}$\square$

\bigskip

With respect to the presentation of $H^{\ast}(G;\mathbb{F}_{p})$ in Lemma 2.3,
the Hopf algebra structure on $H^{\ast}(G;\mathbb{F}_{p})$ has been
determined, see \cite[Remark 4.6]{[DZ3]}\textbf{. }For simplicity we reserve
$x_{i}$ for $r_{p}(x_{i})$.

\bigskip

\noindent\textbf{Lemma 3.3.} \textsl{Let }$(G,p)$\textsl{\ be a pair with }%
$G$\textsl{\ exceptional and }$H^{\ast}(G)$\textsl{\ containing non--trivial
}$p$\textsl{--torsion.}

\begin{enumerate}
\item[(3.4)] \textsl{The algebra} $H^{\ast}(G;\mathbb{F}_{2})$ \textsl{has the
presentation}

$H^{\ast}(G_{2};\mathbb{F}_{2})=\mathbb{F}_{2}[x_{6}]/\left\langle x_{6}%
^{2}\right\rangle \otimes\Delta_{\mathbb{F}_{2}}(\zeta_{3})\otimes
\Lambda_{\mathbb{F}_{2}}(\zeta_{5})$;

$H^{\ast}(F_{4};\mathbb{F}_{2})=\mathbb{F}_{2}[x_{6}]/\left\langle x_{6}%
^{2}\right\rangle \otimes\Delta_{\mathbb{F}_{2}}(\zeta_{3})\otimes
\Lambda_{\mathbb{F}_{2}}(\zeta_{5},\zeta_{15},\zeta_{23})$;

$H^{\ast}(E_{6};\mathbb{F}_{2})=\mathbb{F}_{2}[x_{6}]/\left\langle x_{6}%
^{2}\right\rangle \otimes\Delta_{\mathbb{F}_{2}}(\zeta_{3})\otimes
\Lambda_{\mathbb{F}_{2}}(\zeta_{5},\zeta_{9},\zeta_{15},\zeta_{17},\zeta
_{23})$;

$H^{\ast}(E_{7};\mathbb{F}_{2})=\frac{\mathbb{F}_{2}[x_{6},x_{10},x_{18}%
]}{\left\langle x_{6}^{2},x_{10}^{2},x_{18}^{2}\right\rangle }\otimes
\Delta_{\mathbb{F}_{2}}(\zeta_{3},\zeta_{5},\zeta_{9})\otimes\Lambda
_{\mathbb{F}_{2}}(\zeta_{15},\zeta_{17},\zeta_{23},\zeta_{27})$;

$H^{\ast}(E_{8};\mathbb{F}_{2})=\frac{\mathbb{F}_{2}[x_{6},x_{10}%
,x_{18},x_{30}]}{\left\langle x_{6}^{8},x_{10}^{4},x_{18}^{2},x_{30}%
^{2}\right\rangle }\otimes\Delta_{\mathbb{F}_{2}}(\zeta_{3},\zeta_{5}%
,\zeta_{9},\zeta_{15},\zeta_{23})\otimes\Lambda_{\mathbb{F}_{2}}(\zeta
_{17},\zeta_{27},\zeta_{29})$,
\end{enumerate}

\noindent\textsl{where }

\begin{quote}
$\zeta_{3}^{2}=x_{6}$\textsl{ in }$G_{2},F_{4},E_{6},E_{7},E_{8}%
$\textsl{;\quad}

$\zeta_{5}^{2}=x_{10},\quad\zeta_{9}^{2}=x_{18}$\textsl{ in }$E_{7},E_{8}%
$\textsl{; }

$\zeta_{15}^{2}=x_{30}$, $\zeta_{23}^{2}=x_{6}^{6}x_{10}$ \textsl{in} $E_{8}$,
\end{quote}

\noindent\textsl{on which} \textsl{the action of }$\psi_{2}$ \textsl{is given
by}

\begin{quote}
$\psi_{2}(\zeta_{i})=0$ \textsl{if }$i=3,5,9,17$\textsl{, or }$15,23$\textsl{
for} $F_{4}$;

$\psi_{2}(\zeta_{i})=x_{6}\otimes\zeta_{i-6}$ \textsl{if }$i=15,23$\textsl{
for} $E_{6}$;

$\psi_{2}(\zeta_{15})=\beta_{2}(\zeta_{9}\otimes\zeta_{5})$ \textsl{for}
$E_{7}$;

$\psi_{2}(\zeta_{23})=\beta_{2}(\zeta_{17}\otimes\zeta_{5})$ \textsl{for}
$E_{7}$;

$\psi_{2}(\zeta_{27})=\beta_{2}(\zeta_{17}\otimes\zeta_{9})$ \textsl{for}
$E_{7}$;

$\psi_{2}(\zeta_{15})=\beta_{2}(\zeta_{9}\otimes\zeta_{5})+x_{6}^{2}%
\otimes\zeta_{3}$ \textsl{for} $E_{8}$;

$\psi_{2}(\zeta_{23})=\beta_{2}(\zeta_{17}\otimes\zeta_{5})+%
{\textstyle\sum\limits_{s+t=2}}
x_{6}^{s}\otimes x_{6}^{t}\beta_{2}(\zeta_{5}\otimes\zeta_{5})+x_{10}%
^{2}\otimes\zeta_{3}$ \textsl{for} $E_{8}$;

$\psi_{2}(\zeta_{27})=\beta_{2}(\zeta_{17}\otimes\zeta_{9})+x_{6}^{4}%
\otimes\zeta_{3}$ \textsl{for} $E_{8}$;

$\psi_{2}(\zeta_{29})=x_{10}^{2}\otimes\zeta_{9}+\zeta_{17}\otimes x_{6}%
^{2}+x_{6}^{4}\otimes\zeta_{5}$ \textsl{for} $E_{8}$.
\end{quote}

\begin{enumerate}
\item[(3.5)] \textsl{The algebra} $H^{\ast}(G;\mathbb{F}_{3})$ \textsl{has the
presentation}
\end{enumerate}

\begin{quote}
$H^{\ast}(G_{2};\mathbb{F}_{3})=\Lambda_{\mathbb{F}_{3}}(\zeta_{3},\zeta
_{11})$;

$H^{\ast}(F_{4};\mathbb{F}_{3})=\mathbb{F}_{3}[x_{8}]/\left\langle x_{8}%
^{3}\right\rangle \otimes\Lambda_{\mathbb{F}_{3}}(\zeta_{3},\zeta_{7}%
,\zeta_{11},\zeta_{15})$;

$H^{\ast}(E_{6};\mathbb{F}_{3})=\mathbb{F}_{3}[x_{8}]/\left\langle x_{8}%
^{3}\right\rangle \otimes\Lambda_{\mathbb{F}_{3}}(\zeta_{3},\zeta_{7}%
,\zeta_{9},\zeta_{11},\zeta_{15},\zeta_{17})$;

$H^{\ast}(E_{7};\mathbb{F}_{3})=\mathbb{F}_{3}[x_{8}]/\left\langle x_{8}%
^{3}\right\rangle \otimes\Lambda_{\mathbb{F}_{3}}(\zeta_{3},\zeta_{7}%
,\zeta_{11},\zeta_{15},\zeta_{19},\zeta_{27},\zeta_{35})$;

$H^{\ast}(E_{8};\mathbb{F}_{3})=\mathbb{F}_{3}[x_{8},x_{20}]/\left\langle
x_{8}^{3},x_{20}^{3}\right\rangle \otimes\Lambda_{\mathbb{F}_{3}}(\zeta
_{3},\zeta_{7},\zeta_{15},\zeta_{19},\zeta_{27},\zeta_{35},\zeta_{39}%
,\zeta_{47})$
\end{quote}

\noindent\textsl{on which} \textsl{the action of }$\psi_{3}$ \textsl{is given
by}

\begin{quote}
$\psi_{3}(\zeta_{i})=0$ \textsl{if} $i=3,7,9,17,19$;

$\psi_{3}(\zeta_{11})=-x_{8}\otimes\zeta_{3}$;

$\psi_{3}(\zeta_{15})=-\beta_{3}(\zeta_{7}\otimes\zeta_{7})$;

$\psi_{3}(\zeta_{27})=-\beta_{3}(\zeta_{7}\otimes\zeta_{19})$ \textsl{for}
$E_{7}$;

$\psi_{3}(\zeta_{35})=\zeta_{27}\otimes x_{8}-x_{8}\otimes\zeta_{27}$
$-x_{8}\otimes x_{8}\zeta_{19}$ \textsl{for} $E_{7}$;

$\psi_{3}(\zeta_{27})=\beta_{3}(\zeta_{19}\otimes\zeta_{7})$ \textsl{for}
$E_{8}$;

$\psi_{3}(\zeta_{35})=\zeta_{27}\otimes x_{8}-x_{8}\otimes\zeta_{27}%
-x_{20}\otimes\zeta_{15}-\beta_{3}(x_{8}\zeta_{19}\otimes\zeta_{7})$
\textsl{for} $E_{8}$;

$\psi_{3}(\zeta_{39})=\beta_{3}(\zeta_{19}\otimes\zeta_{19})$;

$\psi_{3}(\zeta_{47})=x_{20}\otimes\zeta_{27}-\zeta_{39}\otimes x_{8}%
-\beta_{3}(x_{20}\zeta_{19}\otimes\zeta_{7}).$
\end{quote}

\begin{enumerate}
\item[(3.6)] \textsl{The algebra} $H^{\ast}(E_{8};\mathbb{F}_{5})$ \textsl{has
the presentation}
\end{enumerate}

\begin{quote}
$H^{\ast}(E_{8};\mathbb{F}_{5})=\mathbb{F}_{5}[x_{12}]/\left\langle x_{12}%
^{5}\right\rangle \otimes\Lambda(\zeta_{3},\zeta_{11},\zeta_{15},\zeta
_{23},\zeta_{27},\zeta_{35},\zeta_{39},\zeta_{47})$
\end{quote}

\noindent\textsl{on which} \textsl{the action of }$\psi_{5}$ \textsl{is given
by}

\begin{quote}
$\psi_{5}(\zeta_{i})=0$, $i=3,11$,

$\psi_{5}(\zeta_{15})=x_{12}\otimes\zeta_{3}$,

$\psi_{5}(\zeta_{23})=2\beta_{5}(\zeta_{11}\otimes\zeta_{11})$,

$\psi_{5}(\zeta_{27})=-x_{12}\otimes\zeta_{15}+2x_{12}^{2}\otimes\zeta_{3}$,

$\psi_{5}(\zeta_{35})=x_{12}\otimes\zeta_{23}+\beta_{5}(3x_{12}\zeta
_{11}\otimes\zeta_{11}-\zeta_{11}\otimes\zeta_{11}x_{12})$,

$\psi_{5}(\zeta_{39})=3x_{12}\otimes\zeta_{27}+x_{12}^{2}\otimes\zeta
_{15}+2x_{12}^{3}\otimes\zeta_{3}$

$\psi_{5}(\zeta_{47})=x_{12}\otimes\zeta_{35}-2x_{12}^{2}\otimes\zeta
_{23}+\beta_{5}(\zeta_{11}\otimes x_{12}^{2}\zeta_{11}$

$\qquad\qquad+3x_{12}\zeta_{11}\otimes x_{12}\zeta_{11}+3x_{12}^{2}\zeta
_{11}\otimes\zeta_{11})$.$\square$
\end{quote}

\noindent\textbf{Remark 3.4.} We emphasize that the ring $H^{\ast}(G)$ (resp.
the algebra $H^{\ast}(G;\mathbb{F}_{p})$) admits many sets of generators
subject to a given degree constraints, and the actions of $r_{p}$ (resp.
$\psi_{p}$) vary sensitively with respect to different choices of a set of
generators. However, in the context of \cite{[DZ2]} the elements
$\varrho_{2i-1}\in E_{3}^{\ast,1}(G;\mathbb{Z})$ and $\zeta_{2s-1}\in
E_{3}^{\ast,1}(G;\mathbb{F}_{p})$ are coherently constructed from explicit
polynomials in certain Schubert classes on $G/T$. Consequently, with respect
to them the presentations of $r_{p}$ in Lemma 3.2 (resp. $\psi_{p}$ in Lemma
3.3) were derived by computing with these polynomials.

\section{The results}

Historically, the Hopf algebras $H^{\ast}(G;\mathbb{F}_{p})$ for exceptional
Lie groups $G$ have been studied by many authors, notably, by Borel, Araki,
Toda, Kono, Mimura, Shimada, see \cite{[DZ3]} for an account for the history.
However, these results can not be directly used in the place of Lemma 3.3
since they were obtained by various methods, presented by generators with
different origins, using case by case computations depending on $G$ and $p$
and without referring to the integral cohomology $H^{\ast}(G)$.

In comparison, with our\textsl{ }generators $\varrho_{j}$'s and $\zeta_{i}$'s
stemming solely from certain polynomials in the Schubert classes on $G/T$ in
\cite{[DZ2]}, the relationships between $H^{\ast}(G)$ and $H^{\ast
}(G;\mathbb{F}_{p})$ are transparent for all prime $p$, as indicated by Lemma
3.2. For this reason the pull--back formula is directly applicable to
translate the Hopf algebra structure on $H^{\ast}(G;\mathbb{F}_{p})$ to the
near--Hopf ring structure on $H^{\ast}(G)$.

An element $x\in H^{\ast}(G)$ is called \textsl{primitive} if $\psi(x)=0$. Let
$\mathcal{P}(G)$ be the graded $\mathbb{Z}$--module of all primitive elements
in $H^{\ast}(G)$.

In Theorems 1--5 below the near--Hopf rings $H^{\ast}(G)$ for all exceptional
$G$ are presented in terms of the set $\mathcal{G}(G)$ of generators specified
in (2.6). In view of Definition 2.4 for the class $\mathcal{C}_{I}%
\in\mathcal{G}(G)$ we note that

\begin{quote}
i) $\mathcal{C}_{I}=x_{i}$ for $I=\{i\}\subseteq e(G,p)$ a singleton;

ii) $\deg\mathcal{C}_{I}=2(i_{1}+\cdots+i_{k})-k+1$ if $I=(i_{1},\cdots
,i_{k})\subseteq e(G,p)$;

iii) the value of $\psi(\mathcal{C}_{I})$, $I\subseteq e(G,p)$, is determined
by $\psi_{p}(\zeta_{2s-1})$, $s\in e(G,p)$ because of the relation $\delta
_{p}\circ\psi_{p}=\psi\circ\delta_{p}$.
\end{quote}

\bigskip

\noindent\textbf{Theorem 1.} \textsl{With respect to the ring presentation}

\begin{quote}
$H^{\ast}(G_{2})=\Delta_{\mathbb{Z}}(\varrho_{3})\otimes\Lambda_{\mathbb{Z}%
}(\varrho_{11})\oplus\tau_{2}(G_{2})$
\end{quote}

\noindent\textsl{where}

\begin{quote}
$\tau_{2}(G_{2})=F_{2}[x_{6}]^{+}/\left\langle x_{6}^{2}\right\rangle
\otimes\Delta_{\mathbb{F}_{2}}(\varrho_{3})$\textsl{, }
\end{quote}

\noindent\textsl{and where }$\varrho_{3}^{2}=x_{6}$\textsl{,} $x_{6}%
\varrho_{11}=0$\textsl{, the reduced co--product }$\psi$\textsl{ is given by}

\begin{quote}
$\{\varrho_{3},x_{6}\}\subset\mathcal{P}(G_{2})$; $\quad$

$\psi(\varrho_{11})=\delta_{2}(\zeta_{5}\otimes\zeta_{5})$.
\end{quote}

\noindent\textbf{Theorem 2.} \textsl{With respect to the ring presentation}

\begin{quote}
$H^{\ast}(F_{4})=\Delta_{\mathbb{Z}}(\varrho_{3})\otimes\Lambda_{\mathbb{Z}%
}(\varrho_{11},\varrho_{15},\varrho_{23})\oplus\tau_{2}(F_{4})\oplus\tau
_{3}(F_{4})$
\end{quote}

\noindent\textsl{where}

\begin{quote}
$\tau_{2}(F_{4})=F_{2}[x_{6}]^{+}/\left\langle x_{6}^{2}\right\rangle
\otimes\Delta_{\mathbb{F}_{2}}(\varrho_{3})\otimes\Lambda_{\mathbb{F}_{2}%
}(\varrho_{15},\varrho_{23})$\textsl{,}

$\tau_{3}(F_{4})=F_{3}[x_{8}]^{+}/\left\langle x_{8}^{3}\right\rangle
\otimes\Lambda_{\mathbb{F}_{3}}(\varrho_{3},\varrho_{11},\varrho_{15}%
)$\textsl{,}
\end{quote}

\noindent\textsl{and} \textsl{where }$\varrho_{3}^{2}=x_{6}$\textsl{, }
$x_{6}\varrho_{11}=0$\textsl{, }$x_{8}\varrho_{23}=0$\textsl{, the reduced
co--product }$\psi$\textsl{ is given by}

\begin{quote}
$\{\varrho_{3},x_{6},x_{8}\}\subset\mathcal{P}(F_{4})$;

$\psi(\varrho_{11})=\delta_{2}(\zeta_{5}\otimes\zeta_{5})+x_{8}\otimes
\varrho_{3}$;

$\psi(\varrho_{15})=-\delta_{3}(\zeta_{7}\otimes\zeta_{7})$

$\psi(\varrho_{23})=\delta_{3}(\zeta_{7}\otimes\zeta_{7}x_{8}-\zeta_{7}%
x_{8}\otimes\zeta_{7})$.
\end{quote}

\noindent\textbf{Theorem 3.} \textsl{With respect to the ring presentation}

\begin{quote}
$H^{\ast}(E_{6})=\Delta_{\mathbb{Z}}(\varrho_{3})\otimes\Lambda_{\mathbb{Z}%
}(\varrho_{9},\varrho_{11},\varrho_{15},\varrho_{17},\varrho_{23})\oplus
\tau_{2}(E_{6})\oplus\tau_{3}(E_{6})$
\end{quote}

\noindent\textsl{where}

\begin{quote}
$\tau_{2}(E_{6})=F_{2}[x_{6}]^{+}/\left\langle x_{6}^{2}\right\rangle
\otimes\Delta_{\mathbb{F}_{2}}(\varrho_{3})\otimes\Lambda_{\mathbb{F}_{2}%
}(\varrho_{9},\varrho_{15},\varrho_{17},\varrho_{23})$\textsl{,}

$\tau_{3}(E_{6})=F_{3}[x_{8}]^{+}/\left\langle x_{8}^{3}\right\rangle
\otimes\Lambda_{\mathbb{F}_{3}}(\varrho_{3},\varrho_{9},\varrho_{11}%
,\varrho_{15},\varrho_{17})$\textsl{,}
\end{quote}

\noindent\textsl{and where }$\varrho_{3}^{2}=x_{6}$\textsl{, }$x_{6}%
\varrho_{11}=0$\textsl{, }$x_{8}\varrho_{23}=0$\textsl{, the reduced
co--product }$\psi$\textsl{ is given by}

\begin{quote}
$\{\varrho_{3},\varrho_{9},\varrho_{17},x_{6},x_{8}\}\subset\mathcal{P}%
(E_{6})$;

$\psi(\varrho_{11})=\delta_{2}(\zeta_{5}\otimes\zeta_{5})+x_{8}\otimes
\varrho_{3}$;

$\psi(\varrho_{15})=x_{6}\otimes\varrho_{9}-\delta_{3}(\zeta_{7}\otimes
\zeta_{7})$;

$\psi(\varrho_{23})=x_{6}\otimes\varrho_{17}+\delta_{3}(\zeta_{7}x_{8}%
\otimes\zeta_{7}-\zeta_{7}\otimes\zeta_{7}x_{8})$.
\end{quote}

\noindent\textbf{Theorem 4.} \textsl{With respect to the ring presentation}

\begin{quote}
$H^{\ast}(E_{7})=\Delta_{\mathbb{Z}}(\varrho_{3})\otimes\Lambda_{\mathbb{Z}%
}(\varrho_{11},\varrho_{15},\varrho_{19},\varrho_{23},\varrho_{27}%
,\varrho_{35})\oplus\tau_{2}(E_{7})\oplus\tau_{3}(E_{7})$
\end{quote}

\noindent\textsl{where}

\begin{quote}
$\tau_{2}(E_{7})=\frac{\mathbb{F}_{2}[x_{6},x_{10},x_{18},\mathcal{C}_{I}%
]^{+}}{\left\langle x_{6}^{2},x_{10}^{2},x_{18}^{2},\mathcal{D}_{J}%
,\mathcal{R}_{K},\mathcal{S}_{I,J},\mathcal{H}_{t,L}\right\rangle }%
\otimes\Delta_{\mathbb{F}_{2}}(\varrho_{3})\otimes\Lambda_{\mathbb{F}_{2}%
}(\varrho_{15},\varrho_{23},\varrho_{27})$

\textsl{\qquad with }$t\in e(E_{7},2)=\{3,5,9\}$, $I,J,L\subseteq e(E_{7}%
,2)$\textsl{, }$\left\vert I\right\vert ,\left\vert J\right\vert \geq
2$\textsl{;}

$\tau_{3}(E_{7})=\frac{\mathbb{F}_{3}[x_{8}]^{+}}{\left\langle x_{8}%
^{3}\right\rangle }\otimes\Lambda_{\mathbb{F}_{3}}(\varrho_{3},\varrho
_{11},\varrho_{15},\varrho_{19},\varrho_{27},\varrho_{35})$\textsl{,}
\end{quote}

\noindent\textsl{and} \textsl{where }$\varrho_{3}^{2}=x_{6}$\textsl{, }%
$x_{8}\varrho_{23}=0$\textsl{,} \textsl{the reduced co--product }$\psi
$\textsl{ is given by}

\begin{quote}
$\{\varrho_{3},x_{6},x_{8},x_{10},x_{18}\}\subset\mathcal{P}(E_{7});$

$\psi(\varrho_{11})=\delta_{2}(\zeta_{5}\otimes\zeta_{5})+x_{8}\otimes
\varrho_{3};$

$\psi(\varrho_{15})=\delta_{2}(\zeta_{9}\otimes\zeta_{5})+\delta_{3}(\zeta
_{7}\otimes\zeta_{7});$

$\psi(\varrho_{19})=\delta_{2}(\zeta_{9}\otimes\zeta_{9});$

$\psi(\varrho_{23})=\delta_{2}(\zeta_{17}\otimes\zeta_{5})+\delta_{3}%
(\zeta_{7}x_{8}\otimes\zeta_{7}-\zeta_{7}\otimes\zeta_{7}x_{8});$

$\psi(\varrho_{27})=\delta_{2}(\zeta_{17}\otimes\zeta_{9})-\delta_{3}%
(\zeta_{7}\otimes\zeta_{19});$

$\psi(\varrho_{35})=\delta_{2}(\zeta_{17}\otimes\zeta_{17})+x_{8}%
\otimes\varrho_{27}-\varrho_{27}\otimes x_{8}+x_{8}\otimes x_{8}\varrho_{19};$

$\psi_{2}(\zeta_{2i-1})=0$, $i\in e(E_{7},2)$.
\end{quote}

\noindent\textbf{Theorem 5.} \textsl{With respect to the ring presentation}

\begin{quote}
$H^{\ast}(E_{8})=\Delta_{\mathbb{Z}}(\varrho_{3},\varrho_{15},\varrho
_{23})\otimes\Lambda_{\mathbb{Z}}(\varrho_{27},\varrho_{35},\varrho
_{39},\varrho_{47},\varrho_{59})\underset{p=2,3,5}{\oplus}\tau_{p}(E_{8})$
\end{quote}

\noindent\textsl{where}

\begin{quote}
$\tau_{2}(E_{8})=\frac{\mathbb{F}_{2}[x_{6},x_{10},x_{18},x_{30}%
,\mathcal{C}_{I}]^{+}}{\left\langle x_{6}^{8},x_{10}^{4},x_{18}^{2},x_{30}%
^{2},\mathcal{D}_{J},\mathcal{R}_{K},\mathcal{S}_{I,J},\mathcal{H}%
_{t,L}\right\rangle }\otimes\Delta_{\mathbb{F}_{2}}(\varrho_{3},\varrho
_{15},\varrho_{23})\otimes\Lambda_{\mathbb{F}_{2}}(\varrho_{27})$
\end{quote}

\noindent\textsl{with} $t\in e(E_{8},2)=\{3,5,9,15\}$\textsl{,} $K,I$%
\textsl{,}$J,L\subseteq e(E_{8},2)$\textsl{, }$\left\vert I\right\vert
,\left\vert J\right\vert \geq2$\textsl{,} $\left\vert K\right\vert \geq3$;

\begin{quote}
$\tau_{3}(E_{8})=\frac{\mathbb{F}_{3}[x_{8},x_{20},\mathcal{C}_{\{4,10\}}%
]^{+}}{\left\langle x_{8}^{3},x_{20}^{3},x_{8}^{2}x_{20}^{2}\mathcal{C}%
_{\{4,10\}},\mathcal{C}_{\{4,10\}}^{2}\right\rangle }\otimes\Lambda
_{\mathbb{F}_{3}}(\varrho_{3},\varrho_{15},\varrho_{27},\varrho_{35}%
,\varrho_{39},\varrho_{47})$;

$\tau_{5}(E_{8})=\frac{\mathbb{F}_{5}[x_{12}]^{+}}{\left\langle x_{12}%
^{5}\right\rangle }\otimes\Lambda_{\mathbb{F}_{5}}(\varrho_{3},\varrho
_{15},\varrho_{23},\varrho_{27},\varrho_{35},\varrho_{39},\varrho_{47})$,
\end{quote}

\noindent\textsl{and where}

\begin{quote}
$\varrho_{3}^{2}=x_{6}$\textsl{, }$\varrho_{15}^{2}=x_{30}$\textsl{, }%
$\varrho_{23}^{2}=x_{6}^{6}x_{10}$\textsl{,}

$x_{2s}\varrho_{3s-1}=0$ \textsl{for} $s=4,5$\textsl{;}

$x_{8}\varrho_{59}=x_{20}^{2}\mathcal{C}_{\{4,10\}}$\textsl{; }$x_{20}%
\varrho_{23}=x_{8}^{2}\mathcal{C}_{\{4,10\}}$\textsl{;}

$x_{12}\varrho_{59}=0$\textsl{,}
\end{quote}

\noindent\textsl{the reduced co--product }$\psi$\textsl{ is given by}

\begin{quote}
$\{\varrho_{3},x_{6},x_{8},x_{10},x_{12},x_{18},x_{20}\}\subset\mathcal{P}%
(E_{8})$;

$\psi(\varrho_{15})=\delta_{2}(\zeta_{9}\otimes\zeta_{5})+x_{6}^{2}%
\otimes\varrho_{3}-\delta_{3}(\zeta_{7}\otimes\zeta_{7})+x_{12}\otimes
\varrho_{3}$;

$\psi(\varrho_{23})=\delta_{2}(\zeta_{17}\otimes\zeta_{5}+%
{\textstyle\sum\limits_{s+t=2}}
x_{6}^{s}\zeta_{5}\otimes x_{6}^{t}\zeta_{5})+x_{10}^{2}\otimes\varrho_{3}$

$\qquad+\delta_{3}(x_{8}\zeta_{7}\otimes\zeta_{7}-\zeta_{7}\otimes\zeta
_{7}x_{8})-\delta_{5}(\zeta_{11}\otimes\zeta_{11})$;

$\psi(\varrho_{27})=\delta_{2}(\zeta_{17}\otimes\zeta_{9})+\delta_{3}%
(\zeta_{19}\otimes\zeta_{7})-x_{12}\otimes\varrho_{15}+(x_{6}^{4}+2x_{12}%
^{2})\otimes\varrho_{3}${\small ;}

$\psi(\varrho_{35})=\delta_{2}(\zeta_{17}\otimes\zeta_{17})-\varrho
_{27}\otimes x_{8}+x_{8}\otimes\varrho_{27}+x_{20}\otimes\varrho_{15}%
+\delta_{3}(x_{8}\zeta_{19}\otimes\zeta_{7})$

$\qquad+2x_{12}\otimes\rho_{23}+\delta_{5}(x_{12}\zeta_{11}\otimes\zeta
_{11}+3\zeta_{11}\otimes\zeta_{11}x_{12})$;

$\psi(\varrho_{39})=\delta_{2}(%
{\textstyle\sum\limits_{s+t=2}}
x_{10}^{s}\zeta_{9}\otimes x_{10}^{r}\zeta_{9})-\delta_{3}(\zeta_{19}%
\otimes\zeta_{19})+x_{12}\otimes\varrho_{27}$

$\qquad+2x_{12}^{2}\otimes\varrho_{15}-x_{12}^{3}\otimes\varrho_{3}$;

$\psi(\varrho_{47})=\delta_{2}(%
{\textstyle\sum\limits_{s+t=6}}
x_{6}^{s}\zeta_{5}\otimes x_{6}^{r}\zeta_{5})-x_{20}\otimes\varrho
_{27}+\varrho_{39}\otimes x_{8}+\delta_{3}(x_{20}\zeta_{19}\otimes\zeta_{7})$

$\qquad+2x_{12}\otimes\varrho_{35}+x_{12}^{2}\otimes\varrho_{23}+\delta
_{5}(\zeta_{11}\otimes x_{12}^{2}\zeta_{11}+%
{\textstyle\sum\limits_{s+t=2}}
x_{12}^{s}\zeta_{11}\otimes x_{12}^{r}\zeta_{11})$;$\ \ $

$\psi(\varrho_{59})=\delta_{2}(x_{10}^{2}\zeta_{29}\otimes\zeta_{9}%
+x_{30}\zeta_{17}\otimes\zeta_{5}x_{6}+x_{18}\zeta_{29}\otimes\zeta_{5}%
x_{6}+x_{6}^{4}\zeta_{29}\otimes\zeta_{5}$

$\qquad+\zeta_{29}\otimes\zeta_{29}+x_{10}^{2}\zeta_{17}\otimes\zeta_{9}%
x_{6}^{2}+\zeta_{17}\otimes x_{6}^{2}\zeta_{29}+x_{6}^{4}\zeta_{17}%
\otimes\zeta_{5}x_{6}^{2}$

$\qquad+x_{18}\zeta_{17}\otimes\zeta_{5}x_{6}^{4}+x_{6}^{4}x_{10}^{2}%
\otimes\zeta_{5}\zeta_{9}+x_{10}^{2}\otimes\zeta_{9}\zeta_{29}+x_{6}%
^{4}\otimes\zeta_{5}\zeta_{29})$

$\qquad+\delta_{3}(%
{\textstyle\sum\limits_{s+t=1}}
(-x_{20})^{s}\zeta_{19}\otimes x_{20}^{r}\zeta_{19})+2\delta_{5}(%
{\textstyle\sum\limits_{s+t=4}}
(-x_{12})^{s}\zeta_{11}\otimes x_{12}^{r}\zeta_{11})$;

$\psi_{p}(\zeta_{2i-1})=0$ \textsl{for}
$(p,i)=(2,3),(2,5),(2,9),(3,4),(3,10),(5,6)$;

$\psi_{2}(\zeta_{29})=x_{10}^{2}\otimes\zeta_{9}+\zeta_{17}\otimes x_{6}%
^{2}+x_{6}^{4}\otimes\zeta_{5}$.
\end{quote}

Explanations on the notation used in Theorems 1--5 are in order.

i) Given a ring $\mathcal{A}$ and a finite set $\{t_{i}\}_{1\leq i\leq n}$ the symbols

\begin{quote}
$\Delta_{\mathbb{Z}}(t_{i})_{1\leq i\leq n}$, $\Lambda_{\mathbb{Z}}%
(t_{i})_{1\leq i\leq n}$ and $\mathcal{A}[t_{i}]_{1\leq i\leq n}^{+}$
\end{quote}

\noindent stand, respectively, for the $\mathbb{Z}$--module with basis the
square free monomials in $t_{1},\cdots,t_{n}$, the exterior algebra generated
by $t_{1},\cdots,t_{n}$ over $\mathbb{Z}$ and the polynomials ring in
$t_{1},\cdots,t_{n}$ over $\mathcal{A}$ of positive degrees in $t_{1}%
,\cdots,t_{n}$.

ii) The relations of the types $\mathcal{D}_{J},\mathcal{R}_{K}$,
$\mathcal{S}_{I,J}$ $\mathcal{H}_{t,I}$ (that occur only in the presentations
of $\tau_{2}(E_{7})$ and $\tau_{2}(E_{8})$) are many. Their full expressions
can be found in \cite[Theorem 6]{[DZ2]}.

Finally, a few words concerning the proofs of Theorems 1--5. Firstly, the
presentation of $H^{\ast}(G)$ as a ring comes from \cite[Theorem 6]{[DZ2]}.
Next, granted with Lemmas 3.2 and 3.3 and taking into account for the obvious relation

\begin{quote}
$\delta_{p}=r_{p}^{-1}\circ\beta_{p}:$ $H^{\ast}(G\times G;\mathbb{F}%
_{p})\rightarrow\tau_{p}(G\times G)$,
\end{quote}

\noindent it is straightforward to apply (3.1) to deduce the expressions for
$\psi(\varrho_{i})$.

\end{document}